\input amstex
\documentstyle{amsppt}
\hsize 6.6in
\loadbold
\hoffset=1cm
\voffset=2cm
\topmatter

\title Computational Geometry in $Heis^3$ 
\endtitle

\author Andrey Marenich\endauthor

\thanks{The author is grateful to Professor Boris Apanasov for supervising this research.
Without his help and encouragement this work would never be possible}\endthanks
\address{Dept of Mathematics, University of Oklahoma, Norman, OK 73019}\endaddress
\keywords Left invariant metric, geodesic lines, Heisenberg geometry \endkeywords
\subjclass  53C15, 53C20 \endsubjclass
\email amarenich\@math.ou.edu\endemail

\abstract Among eight possible geometric structures on
three-dimensional manifolds less studied from the differential
geometric point of view are those modelled on the Heisenberg group
$Heis^3$. We consider the Heisenberg left-invariant metric and use
some results on Levi-Civita connection and curvature tensor to
present solutions of equations for geodesic lines in the
Heisenberg group. Using "Mathematica" software package we also
present drawings of geodesic lines and metric balls in the
Heisenberg group.
\endabstract

\dedicatory \enddedicatory

\endtopmatter

\document
\baselineskip 12pt

\head \bf 1. Heisenberg group $Heis^{2n+1}$ \endhead

Heisenberg group plays an important role in many branches of
mathematics such as representation theory, harmonic analysis, PDEs
or even quantum mechanics, where it was initially defined as a
group of $3\times 3$ matrices
$$
\{m(x,y,t)\}=\left\{ \left(
\matrix 1 & y &  t \\
0 & 1 & x\\
0 & 0 & 1\\
\endmatrix  \right):x,y,t\in \Bbb R\right\}.
$$
with the usual multiplication rule.

We  will use  the following complex definition of the Heisenberg
group:
$$
Heis^{2n+1}= \Bbb C^{n}\times\Bbb R = \{(z,t):z\in \Bbb C,
t\in\Bbb R\} \quad\text{with}\quad (z,t) \cdot (z',t')= (z+z',
t+t'+ Im(<z,z'>))    \tag 1
$$
where $<,>$ is the usual Hermitian product in $\Bbb C^{n}$. See
\cite{A}, \cite{AX}, \cite{G}, \cite{Ma} for details and other
normalizations.

The element zero {\bf 0}=(0,...,0) is the unit of this group
structure and the inverse element for $(z,t)$ is
$(z,t)^{-1}=(-z,-t)$. Let $a=(z,t)$, $b=(w,s)$ and $c=(z',t')$.
The commutator of the elements $a,b\in Heis^{3}$ is equal to
$$
[a,b]=aba^{-1}b^{-1}=(z,t)\cdot(w,s)\cdot(-z,-t)\cdot(-w,-s)=
(z+w-z-w,t+s-t-s+\alpha)=(0,\alpha)
$$
where $\alpha\neq 0$ in general. For example
$[(1,0),(i,0)]=(0,2)\neq(0,0)$. Which shows that $Heis^{3}$ is not
abelian. On the other hand , for any $a,b,c\in Heis^{3}$, their
double commutator is
$$
[[a,b],c]=[(0,\alpha),(z',t')]=(z'-z',\alpha+t'-\alpha-t)=(0,0)
$$
This implies that $Heis^{2n+1}$ is a {\it nilpotent Lie group}
with nilpotency 2.
\newpage

\head \bf 2. Metrics on Heisenberg group $Heis^{2n+1}$ \endhead
There are various different metrics we can define on $Heis^{2n+1}$.

Let us consider (n+1)-dimensional complex hyperbolic space
$X=H^{n+1}_{\Bbb C}$. Every point on the ideal boundary
$p\in\partial X$ can be identified with the class of asymptotic
geodesics, issuing from this point, which defines a natural
fibration of the space. The dual fibration consists of
horospheres, "centered" at $p$ and orthogonal to all such
geodesics. A group of isometries of the space, which fix the point
$p$, can be represented using Ivasava decomposition as
$Isom(X)_{p}=K\cdot A\cdot H$,where $K$ is a group of "rotations"
around some geodesic, $A$ - group of "dilations" along geodesics,
and $H$ is a group of all "translations" along the horosphere. It
can be verified, that $H$ is actually isomorphic to $Heis^{2n+1}$.
This way the usual complex hyperbolic metric in the hyperbolic
space $H^{n+1}_{\Bbb C}$ naturally induces a metric on the
Heisenberg group.

Another possible example of metric on $Heis^{2n+1}$ is a known
{\it Cygan's metric}, defined as:
$$
\rho_{c}((z,t),(z',t'))= \left | {\|(z-z')\|}^{4} +
(t-t'+Im<z,z'>)^{2}\right |^{1\over 4}
$$

We can see from the definition, that when we increase the radius $\lambda$
of a ball in this metric, the ball grows linearly with $\lambda$ "horizontally", and as
$\lambda^{2}$ in the vertical direction. Also, the ball represents a
convex body (in the Euclidean sense).

Another metric on $Heis^{3}$ will be defined below.

\head \bf 3. Left-invariant metric on Heisenberg group $Heis^{3}$ \endhead

Each point $(x,y,t)=(z,t)\in Heis^{3}=\Bbb C\times\Bbb R$ can be
viewed as a translation from {\bf 0} to this point as
$(x,y,t)\cdot(0,0,0)=(x,y,t)$. Then the Euclidean coordinate
directions are translated to $(x,y,t)\cdot(s,0,0)=(x+s,y,t-sy)$,
$(x,y,t)\cdot(0,s,0)=(x,y+s,t+sx)$,
$(x,y,t)\cdot(0,0,s)=(x,y,t+s)$. Differentiating, we obtain the
vector fields
$$
X=(1, 0, -y), \quad Y=(0,1, x), \quad T=(0, 0,1) \tag 2
$$
which are left-invariant vector fields by construction. We define
the left-invariant metric on $Heis^{3}$ by taking $X,Y,T$ as {\bf
the orthonormal frame} in each tangent space
${T_{(x,y,t)}Heis^{3}}$.

\definition{Definition~1} Denote by $g$ the left-invariant metric on $Heis^3$ such that
vector fields $X$, $Y$ and $T$ are orthonormal ones. The corresponding scalar product
we denote as usual by $(\, , \,)$.
\enddefinition

Due to (2) the coordinate vectors are
$$
{{\partial}\over{\partial x}}= (X+y T), \qquad {{\partial}\over{\partial y}}= (Y - xT), \quad
\text{and} \quad  {{\partial}\over{\partial z}}= T.
$$
Since by our choice $\{X,Y,T\}$ is an orthonormal basis, we arrive at the following formula for the metric tensor
of our left-invariant metric in coordinates $(x,y,z)$:
$$
g=  \left(
\matrix 1 + y^2 & -  xy &  y \\
-xy & 1+ x^2 & -x\\
 y  & - x & 1\\
\endmatrix  \right). \tag 3
$$

The following was proved by V. Marenich \cite{Ma}.

\proclaim{Proposition~1} For  the  covariant  derivatives  of  the  Riemannian
connection of the left-invariant metric,  defined  above  the  following  is true:
$$
\nabla=\left(
\matrix
0&T&-Y\\
-T&0& X\\
-Y&X&0\\
\endmatrix
\right),
\tag 4
$$
where the $(i,j)$-element in the table above equals $\nabla_{E_i}E_j$ for our basis
$$
\{E_k, k=1,2,3\}=\{X,Y,T\}.
$$
\endproclaim

\head \bf 4. {Geodesic lines in $Heis^{3}$}\endhead
We can find equations of geodesics issuing from {\bf 0}=(0,0,0) following well know results \cite{Ma}.
Let $c(t)$ be such a geodesics
with a natural parameter $t$,  and its vector of velocity given by
$$
\dot c(t)= \alpha(t)X(t) +\beta(t)Y(t) + \gamma(t)T. \tag 5
$$
Then the equation  of a geodesic $\nabla_{\dot c(t)}{\dot c(t)}\equiv 0$ and our
table of covariant derivatives (3) give:
$$
(\alpha'(t)+2\gamma\beta(t)) X(t) + (\beta'(t)-2\gamma\alpha(t)) Y(t) + \gamma'(t)T =0.
$$
Thus we  easily  obtain  the following equations for coordinates of the
vector of velocity of the geodesic $c(t)$ in our left-invariant moving frame:
$$
\left\{
\matrix
\matrix
{\alpha}'(t)+2\gamma\beta(t)=0\\
\beta '(t)-2\gamma\alpha(t)=0
\endmatrix
& \gamma '(t)=0
\endmatrix
\right.
\tag 6
$$
or
$$
\left\{
\matrix
\matrix
(\alpha(t)+\beta(t))'-2\gamma(\alpha(t)-\beta(t))=0\\
(\alpha(t)-\beta(t))'+2\gamma(\alpha(t)+\beta(t))=0\\
\endmatrix
& \gamma '(t)=0
\endmatrix
\right.
\tag 7
$$
Because the parameter $t$ is natural we have
$$
\alpha^2(t)+\beta^2(t))+\gamma^2\equiv 1,
$$
and we  could  take $\gamma(t)\equiv \gamma$ where $|\gamma|\leq
1$ is the $\cos$ of the angle between $\dot c(0)$ and $T$-axis.
For $|\gamma|=1$ we have "vertical" geodesic, coinciding with
"$z$-axis" in $Heis^{3}$, which is an integral line of the
left-invariant vector field $T$. For $\gamma=0$ our equations are
linear. For $\gamma\not=0$ after some easy computation one could
find that:
$$
\left\{
\matrix
\alpha(t)=r \cos(2\gamma t +\phi)\\
\beta(t)= r \sin(2\gamma t +\phi)
\endmatrix
\right.
\tag 8
$$
where $r=\sqrt{\alpha^2+\beta^2}$.   To   find 
equations  for geodesics $c(t)=(x(t), y(t), z(t))$
issuing from $\text{\bf 0}$ we note that if
$$
\dot c(t)=\alpha(t)X(t) + \beta(t)Y(t) + \gamma(t)T
$$
and our left-invariant vector fields are
$$
X=(1,0, -y), \quad Y=(0,1, x), \quad T=(0, 0,1),
$$
then
$$
{\partial\over\partial x}=X + yT, \quad \text{ and } \quad
{\partial\over\partial y}=Y - xT\,\,.
$$
Therefore we easily have:
$$
\left\{
\matrix
\dot x(t)= \alpha(t)\\
\dot y(t)= \beta(t)   \\
\dot z(t) = \gamma - \alpha(t)y(t) + \beta(t)x(t)\\
\endmatrix
\right. \tag 9
$$
After  some  computations  this gives  the following equations for geodesics issuing from zero (see \cite{Ma}):

\proclaim{Proposition~2}   Geodesic  lines  issuing  from  zero  $\text{\bf 0}$ in the
Heisenberg group $Heis^3$ satisfy to the following equations:
$$
\left\{
\matrix
x(t)= {r\over 2\gamma} (\sin(2\gamma t+\phi)-\sin(\phi))\\
y(t)= {r\over 2\gamma} (\cos(\phi)-\cos(2\gamma t+\phi))\\
z(t)={{1+\gamma^2}\over{2\gamma}}t-
{{1-\gamma^2}\over{4\gamma^2}}\sin(2\gamma t)
\endmatrix
\right.
\tag 10
$$
for some  numbers  $\phi$, $r$  which  could  be  defined from the initial condition
$\dot c(0)=(r\cos(\phi), r\sin(\phi), \gamma)$; or if $\gamma=0$,  then they are "horizontal"
and  satisfy  the  following equations:
$$
\left\{
\matrix
x(t)=\alpha(0)t\\
y(t)=\beta(0)t\\
z(t) = 0\\
\endmatrix
\right.
\tag 11
$$
\endproclaim

To find equations of geodesics issuing from an arbitrary point
$(x,y,z)\in Heis^{3}$ it is sufficient to  use left translation to
this point and apply to the equation above the multiplication rule
(1).

\head \bf 5. Computer generated pictures of geodesic lines and metric balls in $H^{3}$ \endhead

Here we produce several images of geodesics and metric spheres in
the Heisenberg group with the left invariant metric. They are
surprisingly different from their analogs in the metric induced in
the Heisenberg group by the metric of negative sectional curvature
in the complex hyperbolic 2-space and Cygan's metric, see
\cite{AX} and \cite{G}. The reader may compare our images with
many computer generated images in the last geometry given in
Goldman's book \cite{G}.

In Figure 1, we present the exp-image of the $\{X,T\}$-coordinate
plane.


In Figure 2, we present metric balls with the center at {\bf 0} and radii $1$ and $3$.
This shows already a big difference with the ball shapes in the Cygan's
metric, see Section 2.


In Figure 3, we present the half of the ball of radius $5$. 


In Figure 4, we present the amplified singular point of the sphere of radius $5$ and
the amplified neighborhood of this point of the same sphere.


In Figure 5, we present the singular point of the metric sphere of radius $t=20$.

\def\ref#1{[#1]}
\eightpoint
\parindent=36pt

\head  REFERENCES \endhead
\bigskip

\frenchspacing

\item{\ref{A}}
Boris Apanasov, Conformal Geometry of Discrete Groups and Manifolds -- De
Gruyter Exp. in Math., {\bf 32}, Walter de Gruyter, Berlin--New York,
2000, XIV + 523 pp.

\item{\ref{AX}} Boris Apanasov and Xiangdong Xie, Geometrically finite complex hyperbolic manifolds. -
Intern. J. of Math., {\bf 8:6} (1997), 703-757.

\item{\ref{G}} William Goldman,
Complex hyperbolic geometry. - Oxford Math. Monographs, Clarenton Press,
Oxford, 1999.

\item{\ref{Ma}} Valery Marenich, Geodesics in Heisenberg groups. --
Geom. Dedicata {\bf 66} (1997), 175--185.

\item{\ref{Mi}} John Milnor, Curvatures of Left Invariant Metrics on Lie
Groups. -- Advances in Math. {\bf 21} (1976), 293--329.

\item{\ref{S}} Peter Scott, The geometry of 3-manifolds. -- Bull. London
Math. Soc. {\bf 15 } (1983), 401--487.

\enddocument
\bye